\documentclass[smallextended]{article}
\usepackage{amsmath}
\usepackage{amsthm}
\usepackage{amsfonts}
\usepackage{amssymb}
\usepackage{graphicx}
\usepackage{dsfont}
\usepackage[T1]{fontenc}
\usepackage{psfrag}
\usepackage{theoremref}
\usepackage{varioref}
\usepackage[justification=centering]{caption}
\usepackage{caption}
\usepackage{subcaption}
\newtheorem{theorem}{Theorem}[section]

\newtheorem{definition}[theorem]{Definition}

\newtheorem{lemma}[theorem]{Lemma}

\newtheorem{proposition}[theorem]{Proposition}
\newtheorem{remark}{Remark}

\author{ Ahmed El kharroubi and \\Soukaina El masmari}
\title{ON A LOSS STORAGE NETWORK WITH FINITE CAPACITY }
\begin{document}
\maketitle
 \begin{abstract}
 Our goal in this paper is to investigate the fluid picture associated with an open large scale storage network of non reliable file servers with  finite capacity. In this storage system new files can be added and a file with only one copy can be lost or duplicated. The  Skorokhod problem with oblique reflection in a bounded convex domain is used to identify the  fluid limits. Such analysis involves the study of three different regimes, the under-loaded, the critical and the overloaded regime. The overloaded regime is of particular importance. To identify the fluid limits, new martingales are obtained and an averaging property is established. This paper is the continuation of the recent paper \cite{EL}.
\end{abstract}


   \section{Introduction}
   In this paper, we are concerned with an open large-scale storage system with non-reliable file servers in a communication network. The overall storage capacity  is assumed to be limited.

   In the network considered, servers can break down randomly  and  when the disk of a given server breaks down, its files are lost, but can be retrieved on the other servers if copies are available. In order to ensure persistence, a duplication mechanism of files to other servers is then performed. The goal is for each file to have at least one copy available on one of the servers as long as possible. Furthermore  in order to use the bandwidth in an optimal way, there should not be too many copies of a given file so that the network can accommodate a large number of distinct files.
    
  In the system considered here, if there is enough storage capacity, a file with one copy can be duplicated on the other servers aiming to guarantee persistence in the system and new files can be admitted to the system for storage, each with  two copies,  otherwise, if capacity doesn't allow the new files are rejected and the duplication is blocked.

The natural critical parameters of the network are $(N,\mu_N,\lambda_N,\xi_N,F_N)$ where $N$ is  the number of servers, $\mu_N$ is the failure rates of servers, $\lambda_N$ the bandwidth
allocated to files duplication, $\xi_N$ is the bandwidth allocated  to new files admission and $F_N$ the total storage capacity.  In this paper it will be  assumed that the total capacity $F_N$ is proportional to $N$, that is 
  \begin{equation}\label{cond1}
 \lim_{N\rightarrow +\infty} \frac{F_{N}}{N} = \bar{\beta}
\end{equation} 
$\bar{\beta}$ is the average storage capacity per server,   
 and   that the parameters $\xi_N,\; \mu_N,\lambda_N$ are given by 
 $$\lambda_N = \lambda N,\quad \mu_N = \mu,\quad and \quad \xi_N = \xi N$$
  for some positive real constants $\lambda, \xi$ and $\mu$.

 The evolution in time of the number of files having one copy and files having two copies is modeled by two sequences of  stochastic processes which are solutions of some  stochastic differential equations with reflecting boundary. In order to study the qualitative behaviour of the system, these stochastic processes are renormalized by a scaling parameter  $N$. The resulting renormalized processes are the unique solution of a Skorokhod problem involving a sequence of  random measures induced by the process describing the free capacity. Our main result shows that, as the scaling parameter goes to infinity, the sequence of renormalized processes  is relatively compact in the space of $\mathbb{R}^2$-valued right continuous functions  on $\mathbb{R}_+$ with left limits and the limit  of any convergent subsequence is the unique solution of a given deterministic  dynamical system  with reflections at the boundary of a bounded convex subset of $\mathbb{R}^2$ (Theorem~\ref{th1}). Without reflections at the boundary this dynamical system  admits a unique equilibrium point. According to the position of this equilibrium point,  three possible regimes  can therefore be derived : the under-loaded, the overloaded, and the critical regime.
 
  In the under-loaded regime,  the probability of saturation of the system is small, and one can suppose that the capacity of the system is infinite and in this case the fluid limits are explicitly identified in \cite{EL}. 
   
   In the overloaded regime,  the capacity $F_N$ is reached in a finite time. In order to identify the fluid limits, exponential martingales are constructed which are useful in studying the limiting hitting time. Furthermore the analysis involves a stochastic averaging principle with an underlying ergodic Markov process.
 
 In the critical regime, a probabilistic study of fluctuations of the processes around the equilibrium point gives the convergence to a reflected diffusion.
   
 Large -scale storage networks of non-reliable file servers with duplication mechanism have been studied in many papers, see for example \cite{Ramabhadran} and \cite{piconia} and \cite{piconib} where the impact of different replicating functionalities in a distributed system on its reliability is investigated using a simple Markov chain model. The present paper is one of the research articles on the stochastic analysis of unreliable storage systems with duplication mechanisms. The series of articles on this type of research began with the fundamental paper \cite{robert}, in which the authors investigated the evolution of a closed loss storage system and used different time-scales to provide an asymptotic description of the network's decay. 

Within the same context, a recent paper \cite{EL} investigated the storage system of non-reliable file servers with the duplication policy as an open network due to the newly added transition of admitting new files to the system. The asymptotic behaviour of the system is studied under a fluid level, and the explicit expression of the associated fluid limits is obtained by solving a Skorokhod problem in the  orthant $\mathbb{R}_2^+$.   
Nevertheless, in \cite{EL} capacity of the system is assumed to be infinite. And in order to give a complete description of a storage network with loss, duplication and admitting policies which is of  real use in practice, in this paper,  capacity of the system is assumed to be finite and the asymptotic behaviour of the system is also studied under a fluid level.  The associated fluid limits are solutions of a Skorokhod problem in a given bounded convex domain in $\mathbb{R}_+^2$. Unfortunately the resolution of the obtained Skorokhod problem is more complex due to the introduction of the process describing the free capacity of the system noted $(m^N(t))$.

\paragraph*{\bf{Outline of the paper}}.\quad Section 2 introduces the stochastic model considered and establishes the stochastic evolution equations of the Markov processes investigated . In Section 3  the link between the fluid equations and the Skorokhod problem is established. It is shown in Theorem~\ref{th1} that the sequence of the scaled processes converges in distribution to a deterministic function which is the unique solution of a given Skorokhod problem. The under-loaded regime and the critical regimes are studied in section 4 and 6. In Section 5 the overloaded regime is investigated.

\section{Stochastic Model}\label{sec:2} 
 In this paper we consider a large-scale storage system which consists of $N$  servers in a communication network. Let  $F_{N}$ be the  total number of  files that can be stored in these  servers. It will be assumed that $F_{N}$ is finite.   The file storage system operates as follows : As long as the storage capacity is not exceeded new files can be  admitted  and files with one copy can be duplicated.

  For $i\in\{1,2\}$, $X_i^N (t)$ denotes the number of files with $i$ copies present in the network  at time $t$ and $(X_0(t))$ denotes the number of files lost for  good. Let $(m^N(t))$ be the number of free places in  the network at time $t\geq 0$. The sequence of the processes $(m^N(t))$ is defined on $\bar{\mathbb{N}}=\mathbb{N}\cup \{+\infty\}$ and  is  given by 
 \begin{equation}\label{fc}
 m^N(t)=F_{N}-2X_{2}^N (t)-X_{1}^N (t)
 \end{equation}  
  The file duplication and admitting  policies can be described   as follows :  conditionally on $(X_{1}^{N} (t), X_2^N(t))=(x_1,x_2)$ with $x_1>0$ and $ 2x_2+x_1<F_{N}$, a file with one copy gets an additional copy with rate $\frac{\lambda N}{x_1}$. If  $m^N(t)\geq 2$, new files can be stored with rate $\xi N$.  Copies of files disappear independently at rate $\mu$. If the last replica of a given file is lost before being  repaired, the file is then definitively lost.

   All events are supposed to occur after an exponentially distributed time.  The admitting, failure  and the duplication processes are then independent Poisson processes. The process $X^{N}(t)=(X_{1}^{N} (t),X_{2}^{N} (t))$ is then a Markov process on the state space 
 $$\mathcal{D}^N=\{(x_1,x_2)\in\mathbb{N}^2\;|\; 2x_2+x_1\leq F_N \}$$

  For $(x_1,x_2) \in \mathbb{N}^{2}$ the $\mathcal{Q}$-matrix $Q^{N} = (q^{N}(.,.))$ of  $(X^{N}(t))$ is given by

      \begin{equation}\label{eq1} (x_1,x_2) \longrightarrow (x_1,x_2)+
      \left\{  \begin{array}{ll}
       (0,1) \ \  \xi N  \mathds{1}_{\lbrace x_1 + 2 x_2 < F_{N} -1 \rbrace} \\
       (1,-1)\ \  2 \mu  x_2\\
       (-1,1) \ \  \lambda N   \mathds{1}_{\lbrace x_1> 0,   x_1 + 2 x_2 < F_{N} \rbrace}\\
       (-1,0) \ \ \mu x_1
      \end{array}
              \right.
  \end{equation}

   \subsection{Stochastic differential equations}

    The evolution equations associated to the Markov processes $(X_{0}^{N}(t))$, $(X_{1}^{N}(t))$ and $(X_{2}^{N}(t))$  are given  by:
    \begin{equation}\label{eq1}
   X_{0}^{N}(t) = X_{0}^N(0)+ \overset{+ \infty}{\underset{i = 1}{\sum }} \int_{0}^{t} \mathds{1}_{\lbrace i \leq X_{1}^{N}(u^{-}) \rbrace} \mathcal{N}_{\mu,i}(du).
   \end{equation}
  \begin{eqnarray}\label{eq2}
   X_{1}^{N}(t)& = & X_{1}^{N}(0) - \int_{0}^{t} \mathds{1}_{\lbrace  X_{1}^{N}(u^{-}) > 0,2X_2^N(u^-)+X_1^N(u^-)<F_{N} \rbrace} \mathcal{N}_{\lambda N}(du)\\\nonumber &-& \overset{+ \infty}{\underset{i = 1}{\sum }} \int_{0}^{t} \mathds{1}_{\lbrace i \leq X_{1}^{N}(u^{-}) \rbrace} \mathcal{N}_{\mu,i}(du)\\\nonumber
   &+ &\overset{+ \infty}{\underset{i = 1}{\sum }} \int_{0}^{t} \mathds{1}_{\lbrace i \leq X_{2}^{N}(u^{-})  \rbrace} \mathcal{N}_{2 \mu,i}(du). 
 \end{eqnarray}
  
   \begin{eqnarray}\label{eq3}
    X_{2}^{N}(t) &=&  X_{2}^{N}(0) + \int_{0}^{t}\mathds{1}_{\lbrace 2X_2^N(u^-)+X_1^N(u^-)<F_{N}-1 \rbrace}\mathcal{N}_{\xi N}(du) \\\nonumber 
    &-&\overset{+ \infty}{\underset{i = 1}{\sum }} \int_{0}^{t} \mathds{1}_{\lbrace i \leq X_{2}^{N}(u^{-})  \rbrace} \mathcal{N}_{2 \mu,i}(du)\\\nonumber
    &+&\int_{0}^{t} \mathds{1}_{\lbrace X_{1}^{N}(u^{-}) > 0, 2X_2^N(u^-)+X_1^N(u^-)<F_{N} \rbrace}\mathcal{N}_{\lambda N}(du)
   \end{eqnarray}

   where $(\mathcal{N}_{\alpha,i})$ denotes an i.i.d sequence of Poisson processes with parameter $\alpha$. All the sequences of Poisson processes are assumed to be independent. And $x(u^-)=\lim\limits_{\substack{s\to u \\ s<u}}x(s) $

     The  equations \eqref{eq2} and \eqref{eq3} can be rewritten as 
     
      \begin{align}\label{eq4}
  X_{1}^{N}(t) &= X_{1}^{N}(0)+  M_{1}^{N}(t) -\mu \int_{0}^{t} X_{1}^{N}(u) du +2\mu \int_{0}^{t} X_2^{N}(u) du\\\nonumber
  & -\lambda N \int_{0}^{t} \mathds{1}_{\lbrace  X_{1}^{N}(u^{-}) > 0,2X_2^N(u^-)+X_1^N(u^-)<F_{N} \rbrace} du  \\\label{eq5}
  X_2^N(t)&=X_2^N(0)+ M_{2}^{N}(t)-2\mu \int_{0}^{t} X_2^{N}(u) du \\\nonumber
 &+ \xi N \int_{0}^{t} \mathds{1}_{\{2X_2^N(u^-)+X_1^N(u^-)<F_{N}-1\}} du\\\nonumber
 & +\lambda N \int_{0}^{t}\mathds{1}_{\lbrace  X_{1}^{N}(u^{-}) > 0,2X_2^N(u^-)+X_1^N(u^-)<F_{N} \rbrace}du  \nonumber 
  \end{align}
where $(M_1^N(t))$ and $(M_2^N(t))$ are martingales associated to   Markov processes  $(X_{1}^{N}(t))$ and $(X_{2}^{N}(t))$ ( see\cite{robert2} pp 348)  given by :
 
    \begin{align}\label{eq6}
     M_{1}^{N}(t) &= \overset{+ \infty}{\underset{i = 1}{\sum }} \int_{0}^{t} \mathsf{1}_{\lbrace i \leq X_{2}^{N}(u^{-}) \rbrace} [ \mathcal{N}_{2 \mu,i}(du) -2 \mu du] \\\nonumber
     &- \int_{0}^{t} \mathds{1}_{\{X_1^N(u)>0,2X_2^N(u)+X_1^N(u)<F_{N}\}} [ \mathcal{N}_{\lambda N}(du) - \lambda N du ]\\\nonumber
     &- \overset{+ \infty}{\underset{i = 1}{\sum }} \int_{0}^{t} \mathsf{1}_{\lbrace i \leq X_{1}^{N}(u^{-}) \rbrace} [ \mathcal{N}_{\mu,i}(du) - \mu du] 
     \end{align}
   \begin{align}\label{eq7}
   M_{2}^{N}(t) &= \int_{0}^{t} \mathds{1}_{\lbrace 2X_2^N(u)+X_1^N(u)<F_{N}-1\rbrace}[ \mathcal{N}_{\xi N}(du) - \xi N du ] \\\nonumber
   &+\int_{0}^{t} \mathds{1}_{\{X_1^N(u)>0,2X_2^N(u)+X_1^N(u)<F_{N}\}} [ \mathcal{N}_{\lambda N}(du) - \lambda N du ] -\\\nonumber   
  &-\overset{+ \infty}{\underset{i = 1}{\sum }} \int_{0}^{t} \mathds{1}_{\lbrace i \leq X_{2}^{N}(u^{-}) \rbrace} [ \mathcal{N}_{2 \mu,i}(du) - 2 \mu du] 
   \end{align}

      The predictable  increasing processes associated to the martingales $(M_{1}^{N}(t))$ and $(M_{2}^{N}(t))$ are respectively given by
      \begin{align}\label{eq8}
       \langle M_{1}^{N}\rangle(t) &= 2 \mu \int_{0}^{t}  X_{2}^{N}(u) du + \mu \int_{0}^{t}  X_{1}^{N}(u) du \\ \nonumber
       &+ \lambda N \int_{0}^{t}\mathds{1}_{\lbrace  X_{1}^{N}(u) > 0,\; 2X_2^N(u)+X_1^N(u)<F_{N} \rbrace} du
       \end{align}
        \begin{align}\label{eq9}
       \langle M_{2}^{N}\rangle(t)& = \xi N \int_{0}^{t} \mathds{1}_{\lbrace  2X_2^N(u)+X_1^N(u)<F_{N}-1 \rbrace} du +2 \mu \int_{0}^{t}  X_{2}^{N}(u) du\\\nonumber
        & + \lambda N \int_{0}^{t}\mathds{1}_{\lbrace  X_{1}^{N}(u) > 0,\;2X_2^N(u)+X_1^N(u)<F_{N} \rbrace} du 
        \end{align}
   
 \section{ Fluid equations and Skorokhod Problem}
 
Let $\mathcal{S}$ be the convex domain in $ \mathbb{R}^2$ given by  $$\mathcal{S}=\{(x_1,x_2)\in\mathbb{R}^2 |x_1\geq 0,x_2\geq 0, 2x_2+x_1\leq \bar{\beta}\}$$ and  $\mathcal{D}(\mathbb{R}_+,\mathbb{R}^2)$ the space of $\mathbb{R}^2$-valued right continuous functions  on $\mathbb{R}_+$ with left limits. Let $\mathcal{M}_{m,n}(\mathbb{R})$ be  the space of $m\times n $ matrices over $\mathbb{R}$. 

 In this paper we consider the following Skorokhod problem in the convex domain $\mathcal{S}$.
 Let   $\theta\in\mathcal{M}_{2,1}(\mathbb{R})$,\quad  $A\in\mathcal{M}_{2,2}(\mathbb{R})$ and  $ R\in\mathcal{M}_{2,2}(\mathbb{R})$. Let  $\nu$ be  the  measure on $[0,+\infty[\times \bar{\mathbb{N}}$ satisfying $\nu([0,t]\times \bar{\mathbb{N}})=t$ for all $t\geq 0$.

\begin{definition}\label{SP}
The couple of functions  $z\in\mathcal{D}(\mathbb{R}_+,\mathbb{R}^2)$ \quad and  $y\in\mathcal{D}(\mathbb{R}_+,\mathbb{R}^2)$ with $z(0)\in \mathcal{S}$, is called  the solution of the  Skorokhod problem associated with the  data $(\theta,\nu,A,R,\mathcal{S})$ and the function
  
\begin{equation}\label{skeq1}
x(t)= z(0)+t\theta+\mathcal{V}(t,\Gamma)+\int_{0}^t Az(s)ds 
\end{equation}
where for $\Gamma$ in a $\sigma$-algebra $\mathcal{B}(\bar{\mathbb{N}})$ $$\mathcal{V}(t,\Gamma)=\begin{pmatrix} 0\\\nu(t,\Gamma)\end{pmatrix}$$

  if the three following conditions hold :
 \begin{enumerate}
 \item \begin{equation}\label{skeq2}
 z(t)=z(0)+t\theta+\mathcal{V}(t,\Gamma)+\int_{0}^t Az(s)ds+Ry(t)
\end{equation}  
 \item   $z(t)\in \mathcal{S}$ for all $t\geq 0$ 
   \item for $i=1,2$ the component $y_i$\; of the function $y$  are  non-decreasing functions with $y_i(0)=0$, and for $t\geq 0$
   \begin{align}  
    &y_1(t)=\int_{0}^{t}\mathds{1}_{\{z_1(s)=0\}} dy_1(s)\\
    &y_2(t)=\int_{0}^{t}\mathds{1}_{\{z_1>0,z_1(s)+2z_2(s)=\bar{\beta}\}} dy_2(s)
    \end{align}
\end{enumerate}

\end{definition}

If $z\in\mathcal{D}(\mathbb{R}_+,\mathbb{R}^d)$ \quad and  $y\in\mathcal{D}(\mathbb{R}_+,\mathbb{R}^d)$ with $z(0)\in \mathcal{S}$ is a solution of the above Skorokhod Problem  then the function $z=(z(t))$ has the following properties. First $z$ behaves on the interior of the set $S$ like a solution of the following  ordinary differential equation  
\begin{equation}\label{ode}
x(t)=x(0)+\theta t+\mathcal{V}(t,\Gamma)+\int_0^t Ax(s)ds
\end{equation}
And second, $z$ is reflected instantaneously at the boundaries $(\partial\mathcal{S})_1=\{x_1=0\}$ and  $(\partial\mathcal{S})_2=\{x_1+2x_2=\bar{\beta}\}$ of the set $\mathcal{S}$. The direction of the reflection on the boundary $(\partial\mathcal{S})_1$  is the first column vector of the reflection matrix $R$ and the direction of reflection on $(\partial\mathcal{S})_2$ is the second column vector the matrix $R$. See for example \cite{tanaka}.

\subsection{Fluid equations}
If $(X^N(t))$ is a sequence of processes, one defines the renormalized sequence of processes of $(X^N(t))$ by  $$\bar{X}^N(t)\stackrel{def}{=} \frac{X^{N}(t)}{N},\; \text{for}\; t\geq 0$$.

From equations \eqref{fc}, \eqref{eq4},  \eqref{eq5} one gets the fluid stochastic differential equations associated with the sequence of  processes $(\bar{X}_1^N(t))$ and $(\bar{X}_2^N(t))$ 

\begin{equation}\label{eq11}
\begin{split}
  \bar{X}_{1}^{N}(t) = & \bar{X}_{1}^{N}(0)+  \bar{M}_{1}^{N}(t)-\lambda t-\mu \int_{0}^{t}  \bar{X}_{1}^{N}(u) du \\
  &\quad +2\mu \int_{0}^{t}  \bar{X}_2^{N}(u) du  + \lambda \int_{0}^{t} \mathds{1}_{\lbrace  \bar{X}_{1}^{N}(u) > 0, m^N(u)=0\rbrace} du \\
  &\qquad+\lambda \int_{0}^{t} \mathds{1}_{\lbrace  \bar{X}_{1}^{N}(u)=0\rbrace} du
  \end{split}  
  \end{equation}
  \begin{equation}\label{12}
  \begin{split}
   \bar{X}_2^N(t) =& \bar{X}_2^N(0) +   \bar{M}_{2}^{N}(t)+(\lambda+\xi)t-2\mu \int_{0}^{t}  \bar{X}_2^{N}(u) du\\
   &\quad -\xi  \int_{0}^{t} \mathds{1}_{\{m^N(u)\leq 1\}} du 
  -\lambda \int_{0}^{t} \mathds{1}_{\lbrace  \bar{X}_{1}^{N}(u)> 0 ,m^N(u)=0\rbrace} du\\
  &\qquad- \lambda \int_{0}^{t} \mathds{1}_{\lbrace  \bar{X}_{1}^{N}(u)= 0\rbrace} du
   \end{split}
  \end{equation}

  The process $(m^N(t))$ evolves on a very rapid time-scale compared with the process $\bar{X}^N(t)\stackrel{def}{=}( \bar{X}_1^N(t), \bar{X}_2^N(t))$. One can see that, while the velocity of the process $(\overline{X}^N(t))$ is of the order O(1), velocity of the process $(m^N(t))$ is much faster than $(\overline{X}^N(t))$ and is of the order O(N).
  
   We consider as in Hunt and Kurtz \cite{HK} the random measure $\nu^N$ on $[0,+\infty[\times \bar{\mathbb{N}}$ defined by 
  \begin{equation}\label{14}
  \nu^N((0,t)\times \Gamma)=\int_0^t\mathds{1}_{\{m^N(u)\in\Gamma\}}du
  \end{equation}
  for all $t\in [0,+\infty[$ and $\Gamma$ in a $\sigma$-algebra $\mathcal{B}(\bar{\mathbb{N}})$. Note that the measure $\nu^N$ satisfies the condition $\nu^N((0,t)\times \bar{\mathbb{N}})=t$. There is a subsequence of the sequence $(\nu^N)$ that converges in distribution to random measure $\nu$ satisfying  $\nu((0,t)\times ~\bar{\mathbb{N}})=t$. (see \cite{HK} for more details). In terms of the random measure $\nu^N$ equations \eqref{eq11},\eqref{12}  becomes 
  \begin{equation}\label{eq15}
\begin{split}
  \bar{X}_{1}^{N}(t) = & \bar{X}_{1}^{N}(0)+  \bar{M}_{1}^{N}(t)-\lambda t-\mu \int_{0}^{t}  \bar{X}_{1}^{N}(u) du \\
  &\quad +2\mu \int_{0}^{t}  \bar{X}_2^{N}(u) du  + \lambda \int_{0}^{t} \mathds{1}_{\lbrace  \bar{X}_{1}^{N}(u) > 0, m^N(u)=0\rbrace} du \\
  &\qquad+\lambda \int_{0}^{t} \mathds{1}_{\lbrace  \bar{X}_{1}^{N}(u)=0\rbrace} du
  \end{split}  
  \end{equation}
  \begin{equation}\label{eq16}
  \begin{split}
   \bar{X}_2^N(t) =& \bar{X}_2^N(0) +   \bar{M}_{2}^{N}(t)+(\lambda+\xi)t-2\mu \int_{0}^{t}  \bar{X}_2^{N}(u) du\\
   &\quad -\xi  \nu^N([0,t]\times\{0,1\}) 
  -\lambda \int_{0}^{t} \mathds{1}_{\lbrace  \bar{X}_{1}^{N}(u)> 0 ,m^N(u)=0\rbrace} du\\
  &\qquad- \lambda \int_{0}^{t} \mathds{1}_{\lbrace  \bar{X}_{1}^{N}(u)= 0\rbrace} du
   \end{split}
  \end{equation}
  
 The above equations can be rewritten in the matrix form as follows
\begin{equation}\label{17}
\begin{split}
\bar{X}^N(t)&=\bar{X}^N(0)+\bar{M}^N(t)+t\bar{\theta}-\xi\mathcal{V}^N(t,\{0,1\})\\
&\qquad+ \int_{0}^t A\bar{X}^N(s)ds+RY^N(t)
\end{split}
\end{equation}

where 
$$\bar{X}^N(t)=\begin{pmatrix}\bar{X}_{1}^{N}(t)\\ \bar{X}_{2}^N(t)\end{pmatrix}\;  \bar{M}^N(t)=\begin{pmatrix} \frac{M_{1}^{N}(t)}{N}\\ \frac{M_{2}^{N}(t)}{N} \end{pmatrix}$$
$$\bar{\theta}=\begin{pmatrix} -\lambda\\\xi +\lambda\end{pmatrix},
\; A=\begin{pmatrix} -\mu & 2\mu \\ 0 & -2\mu \end{pmatrix},\;\; R=\begin{pmatrix} \lambda & \lambda\\ -\lambda & -\lambda\end{pmatrix}$$ 

$$\mathcal{V}^N(t,\{0,1\})=\begin{pmatrix}  \nonumber  0 \\   \nu^N([0,t] \times \{0,1\})\end{pmatrix}$$
$$ Y^N(t)=\begin{pmatrix}
 \int_{0}^{t} \mathds{1}_{\{\bar{X}_1^N(u)= 0\}}du\\
\int_{0}^{t} \mathds{1}_{\lbrace  \bar{X}_{1}^{N}(u)> 0 ,m^N(u)=0\rbrace} du\end{pmatrix}$$

\begin{figure}[hbtp]
\centering
\includegraphics[scale=0.2]{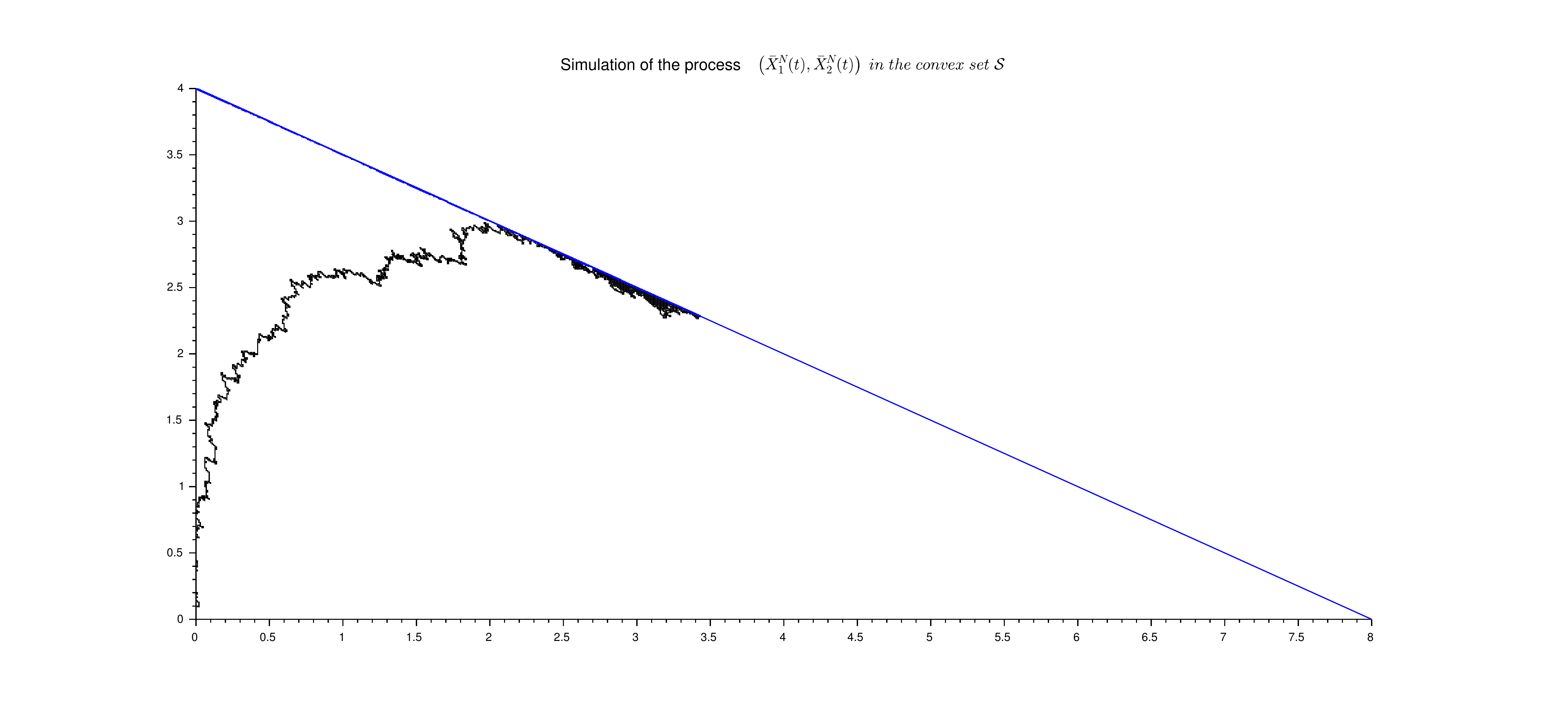}
\caption{}
\label{fig1}
\end{figure}

As illustrated in Figure\ref{fig1} the couple of processes $( \bar{X}^N(t) )$ and $(Y^N (t))$ can be interpreted as the solution of the Skorokhod problem associated with data $(\bar{\theta} ,\nu^N, A,R,\mathcal{S})$ and 

\begin{equation}\label{18}  \bar{V}^N(t)=\bar{X}^N(0)+\bar{M}^N(t)+t\bar{\theta}-\xi\mathcal{V}^N(t,\Gamma)+ \int_{0}^t A\bar{X}^N(s)ds \end{equation}

In the next theorem we prove relative compactness of the sequence of processes 
$\left(\overline{X}^{N}(.),Y^{N}(.), \mathcal{V}^N(.)\right)$   in $\mathcal{D}(\mathbb{R}_+,\mathbb{R}^2)\times \mathcal{M}_{1}(\mathbb{R}_+\times\bar{\mathbb{N}})$. Where $\mathcal{M}_{1}(\mathbb{R}_+\times\bar{\mathbb{N}})$ is the space of Radon measures on $\mathbb{R}_+\times\bar{\mathbb{N}}$.
\begin{theorem}\label{th1} Suppose that 
\begin{eqnarray}\nonumber \underset{N \rightarrow + \infty}{lim} ( \overline{X}_1^N(0), \overline{X}_2^N(0)) = (x_1, x_2) \in \mathcal{S}, 
\end{eqnarray} 
the sequence $\left(\overline{X}^{N}(.),Y^{N}(.), \nu^N(.)\right)$ is then relatively compact in $\mathcal{D}(\mathbb{R}_+,\mathbb{R}^3)$ and the limit $\left( x(.) ,y(.), \nu(.)\right)$ of any convergent subsequence satisfies: 
\begin{align}\label{19}x_1(t)  &=  x_1 - \lambda t-\mu \int_{0}^t x_1(s)ds  + 2\mu \int_{0}^t x_2(s) ds \\\nonumber
   &+ \lambda \int_{[0,t] \times \mathbb{N}} \mathds{1}_{\lbrace x_1(s) > 0 \rbrace} \mathds{1}_{\lbrace  0 \rbrace} (u) \nu (ds \times du)   +\lambda y_1(t)
 \end{align}
\begin{align} \label{20}  x_2(t)  &=  x_2 + (\lambda+\xi) t -2  \mu \int_{0}^t x_2(s)  ds- \xi  \nu ([0,t] \times \{0,1\})  \\\nonumber
   &  - \lambda \int_{[0,t] \times  \mathbb{N}}  \mathds{1}_{\lbrace x_1(s) > 0 \rbrace} \mathds{1}_{\lbrace  0 \rbrace} (u)  \nu (ds \times du) - \lambda y_1(t)  
 \end{align}
 where the function $y_1$  is a non-decreasing function with $y_1(0)=0$, and for $t\geq 0$
   \begin{align}  
    y_1(t)=\int_{0}^{t}\mathds{1}_{\{x_1(s)=0\}} dy_1(s)
    \end{align}
\end{theorem}
 \begin{lemma}\label{lemma}
    
     The  sequences of processes $(\frac{M_{1}^{N}(t)}{N})_{t\geq0}$ and  $(\frac{M_{2}^{N}(t)}{N})_{t\geq0}$ converge in distribution to 0 uniformly on compact sets.  
   \end{lemma}
   \begin{proof}
  Doob's inequalities  show that, for $\epsilon >0$ and $t \geq 0$
     
     $$\mathbb{P}\left( \underset{0 \leq s\leq t }{sup}\frac{M_{i}^{N}(s)}{N} \geq \epsilon \right)\leq  \frac{1}{\epsilon^{2} N^{2}} \mathbb{E} ( \langle M_{i}^{N}\rangle(t))$$
 From equations \eqref{eq8}, \eqref{eq9} one gets 
    \[
    \begin{split}
    \mathbb{E} ( \langle M_{1}^{N}\rangle(t))&\leq \mu F_N+\lambda Nt\\
           \mathbb{E} ( \langle M_{2}^{N}\rangle(t))&\leq \mu F_N+(\lambda+\xi) Nt\\
    \end{split}
\]      
Then from \eqref{cond1} the  sequences of processes $(\frac{M_{1}^{N}(t)}{N})_{t\geq0}$ and  $(\frac{M_{2}^{N}(t)}{N})_{t\geq0}$ converge in distribution to 0 uniformly on any bounded time interval.
\end{proof}
\begin{proof}[ proof of Theorem~\ref{th1}]
 First we prove the relative compactness of the process $$\bar{X}^N(t)=\begin{pmatrix}\bar{X}_{1}^{N}(t)\\ \bar{X}_{2}^N(t)\end{pmatrix}$$  For this we prove separately  that$(\overline{X}_1^{N}(t))$ and $(\overline{X}_2^{N}(t))$ are tight.
  
  For $T>0, \delta>0$ we denote by $\omega_{g}^T(\delta)$ the modulus of continuity of the function g  on $[0,T]$ :
     \begin{equation}
     \omega_{g}^T(\delta)=\sup_{0\leq s\leq t\leq T,|t-s|\leq \delta}\vert g(t)-g(s)\vert
     \end{equation}
 The equation \eqref{eq15} shows that the processes $(\overline{X}_1^{N}(t),Y_1^N(t))$ with $Y_1^N(t)=\lambda\int_{0}^{t} \mathds{1}_{\{X_1^N(u)= 0\}} du)$ is  the unique  solution of  the Skorokhod problem associated to the process 
 \begin{align} \label{eq21}\overline{V}_1^N(t)& =   \overline{X}_1^N(0)  + \overline{M}_1^N(t) - \lambda t + \mu \int_{0}^t (2  \overline{X}_2^N(u)  -  \overline{X}_1^N(u) ) du \\\nonumber
   &+ \lambda \int_{[0,t] \times \mathbb{N}} \mathds{1}_{\lbrace {X_1}^{N}(u) > 0 \rbrace}\mathds{1}_{\lbrace  0 \rbrace} (y)  \nu^N (ds \times dy)   
 \end{align}
 
  By using explicit representation of the solution of the Skorokhod in dimension 1 \cite{ELkaroui}, one has
\[ \|\overline{X}_1^N\|_{\infty,t}\stackrel{def}{=}\sup_{0\leq s\leq t}|\overline{X}_1^N(s)|\leq 2 \|\overline{V}_1^{N}\|_{\infty,t} \] 
and 
\[|\lambda Y_{1}^N(t)|\leq \|\overline{V}_1^{N}\|_{\infty,t}\]
  
By equation (\ref{eq21}), one gets that \begin{align}\nonumber
\|\overline{V}_1^{N}\|_{\infty,t}&\leq  |\overline{X}_1^N(0))|+ 2 \lambda t+\mu \int_0^t\|\overline{X}_1^N\|_{\infty,s}ds\\\nonumber
&+2\mu\int_0^t\|\overline{X}_2^N\|_{\infty,s}ds+ \|\overline{M}_1^{N}\|_{\infty,t} 
\end{align}
and using inequalities given above,
\begin{align}\nonumber
\|\overline{X}_1^N\|_{\infty,t}&\leq 2 |\overline{X}_1^N(0)|+4 \lambda t+2\mu\int_0^t\|\overline{X}_1^N\|_{\infty,s}ds\\\nonumber
&+4\mu\int_0^t\|\overline{X}_2^N\|_{\infty,s}ds+ 2\|\overline{M}_1^{N}\|_{\infty,t}\nonumber
\end{align}
 
\begin{align}\nonumber
\|\overline{X}_2^N\|_{\infty,t}&\leq |\overline{X}_1^N(0)|+|\overline{X}_2^N(0)|+(4\lambda+3\xi) t
+ \|\overline{M}_1^{N}\|_{\infty,t}+\|\overline{M}_2^{N}\|_{\infty,t}\\\nonumber
&+\mu\int_0^t\|\overline{X}_1^N\|_{\infty,s}ds+4\mu\int_0^t\|\overline{X}_2^N\|_{\infty,s}ds
\end{align}
Then 
\begin{equation}\nonumber
\|\overline{X}_1^N\|_{\infty,t}+\|\overline{X}_2^N\|_{\infty,t}\leq H^N(T)+8\mu\int_0^t(\|\overline{X}_1^N\|_{\infty,s}+\|\overline{X}_2^N\|_{\infty,s})ds
\end{equation} 
with \[H^N(t)=3|\overline{X}_1^N(0)|+ |\overline{X}_2^N(0)|+(8\lambda+3\xi)T+3\|\overline{M}_1^{N}\|_{\infty,T}+\|\overline{M}_2^{N}\|_{\infty,T}\]
Gronwall's lemma gives that the relation 
\[\|\overline{X}_1^N\|_{\infty,t}+\|\overline{X}_2^N\|_{\infty,t}\leq H^N(T)e^{8\mu t}\]
holds for all $t\in[0,T]$.
The convergence of martingales and of $|\overline{X}_1^N(0)|$, $|\overline{X}_2^N(0)|$ shows that the sequence $(H^N(T))$ converges in distribution. Consequently for $\epsilon >0$, there exists some $C>0$ such that for $i=1,2$ and all $N\in\mathbb{N}$
\[\mathbb{P}\left(\|\overline{X}_1^N\|_{\infty,t} + \|\overline{X}_2^N\|_{\infty,t}>C\right)\leq \epsilon.
\]

If $\eta>0$, there exists $N_1$ and $\delta>0$ such that for all $N\geq N_1$
\[\delta(\lambda+4\mu C)\leq \frac{\eta}{2}\]
and  \[\mathbb{P}\left(\omega_{\overline{M}_{1}^{N}}^T(\delta)\geq \frac{\eta}{2}\right)\leq \epsilon\]
One gets finally

\begin{align*}
\mathbb{P}\left(\omega_{\overline{V}_{1}^{N}}^T(\delta)\geq \eta\right)&\leq\mathbb{P}\left(2\lambda\delta+2\mu\delta(\|\overline{X}_1^N\|_{\infty,T}+\|\overline{X}_2^N\|_{\infty,T})\geq \frac{\eta}{2}\right)\\
&+\mathbb{P}\left(\omega_{\overline{M}_{1}^{N}}^T(\delta)\geq \frac{\eta}{2}\right)\leq 3\epsilon
\end{align*}

Consequently the sequence $(\overline{V}_{1}^{N}(t))$ is tight and by continuity of the solution of the Skorokhod problem in dimension 1 the sequences $(\overline{X}_1^N(t))$ and $(\overline{Y}_1^{N}(t))$ are tight \cite{Billingsley}.

From equation (\eqref{eq16}) one gets for $s<t$ : 

\begin{align}\nonumber |\overline{X}_2^N(t)-\overline{X}_2^N(s)| &\leq (\lambda+\xi)(t-s)+2\mu\int_s^t|\overline{X}_2^N(u)|du+|\overline{M}_{2}^{N}(t)-\overline{M}_{2}^{N}(s)| \\\nonumber
&+ (2\lambda+\xi)(t-s)+\lambda (Y_1^N(t)-Y_1^N(s)) 
\end{align} and
\begin{align}\mathbb{P}\left(\omega_{\overline{X}_2^N}^T(\delta) \geq \eta\right) &\leq \mathbb{P}\left(\omega_{\overline{M}_{2}^{N}}^T(\delta)\geq \eta/3\right)
+\mathbb{P}\left(\omega_{Y_{1}^{N}}^T(\delta)\geq \eta/3\right)\\\nonumber
&+\mathbb{P}\left(2\mu\delta\|\overline{X}_2^N\|_{\infty,T}+\delta(2 \lambda+3 \xi) \geq \eta/3\right)
\end{align}
There exists $N_1\geq 0$ such that $\delta(2\mu C-\xi)\leq\epsilon$ and 
 \[\mathbb{P}\left(\omega_{\overline{M}_{2}^{N}}^T(\delta)\geq \frac{\eta}{3}\right)\leq \epsilon\]
 and 
  \[\mathbb{P}\left(\lambda\omega_{\frac{Y_{1}^{N}}{N}}^T(\delta)\geq \frac{\eta}{3}\right)\leq \epsilon\]
  and, consequently 
  \[\mathbb{P}\left(\omega_{\overline{X}_2^N}^T(\delta)\geq \eta\right)\leq 3\epsilon\]
  The sequence $\overline{X}_2^N$ is therefore tight.
  
 It remains to prove relative compactness of the sequence of random measures $\nu^N$ on $[0,+\infty[\times\bar{\mathbb{N}}$ . Since for all $N$ and all $t\geq 0$ 
 $$\nu^N([0,t[\times\bar{\mathbb{N}})=t$$
The result is given in \cite{Kurtz} Lemma 1.3. 
\end{proof}
\begin{remark}
The dynamical system given by \eqref{19} and \eqref{20} without constraints can be written as follows 
 \begin{equation*}
  \left\{\begin{array}{cc}
  x_1(t)=&x_1-\lambda t-\mu\int_0^tx_1(s)ds+2\mu\int_0^tx_2(s)ds\\
  x_2(t)=&x_2+(\lambda+\xi) t-2\mu\int_0^tx_2(s)ds\\
  \end{array}
  \right.
  \end{equation*} 
  The unique solution of this ordinary differential equations   noted  $\overline{x}(t)=(\overline{x}_1(t),\overline{x}_2(t))$ is given by   :
\begin{align}\label{22}
\left\{\begin{array}{ll}
&\overline{x}_1(t)=\left(x_1+2x_2-\dfrac{\lambda+2\xi}{\mu}\right)e^{-\mu t}-\left(2x_2-\dfrac{\lambda+\xi}{\mu}\right)e^{-2\mu t}+\dfrac{\xi}{\mu}\\
&\overline{x}_2(t)=\dfrac{\lambda+\xi}{2\mu}+\left(x_2-\dfrac{\lambda+\xi}{2\mu}\right)e^{-2\mu t}
\end{array}
\right.
\end{align}
This dynamical system  admits  unique equilibrium point  $$(\frac{\xi}{\mu},\frac{\lambda+\xi}{2\mu})$$
Thus, according to the position of this equilibrium point in the convex set $\mathcal{S}$, three possible regimes can be considered. Let $$\rho \stackrel{def}{=}  \frac{\lambda+ 2 \xi}{\mu}$$
 The under-loaded regime ($\rho < \bar{\beta}$), the critical regime ($\rho = \bar{\beta}$) and the overloaded regime ($\rho> \bar{\beta}$). Each of the aforementioned regimes will be developed  in detail in the next sections.

\end{remark}
   
\section{The under-loaded regime}
  Throughout this section, we assume that the condition
\begin{equation}\label{23}
\rho<\bar{\beta}
\end{equation}
holds. 

\begin{figure}[hbtp]
\centering
\includegraphics[scale=0.2]{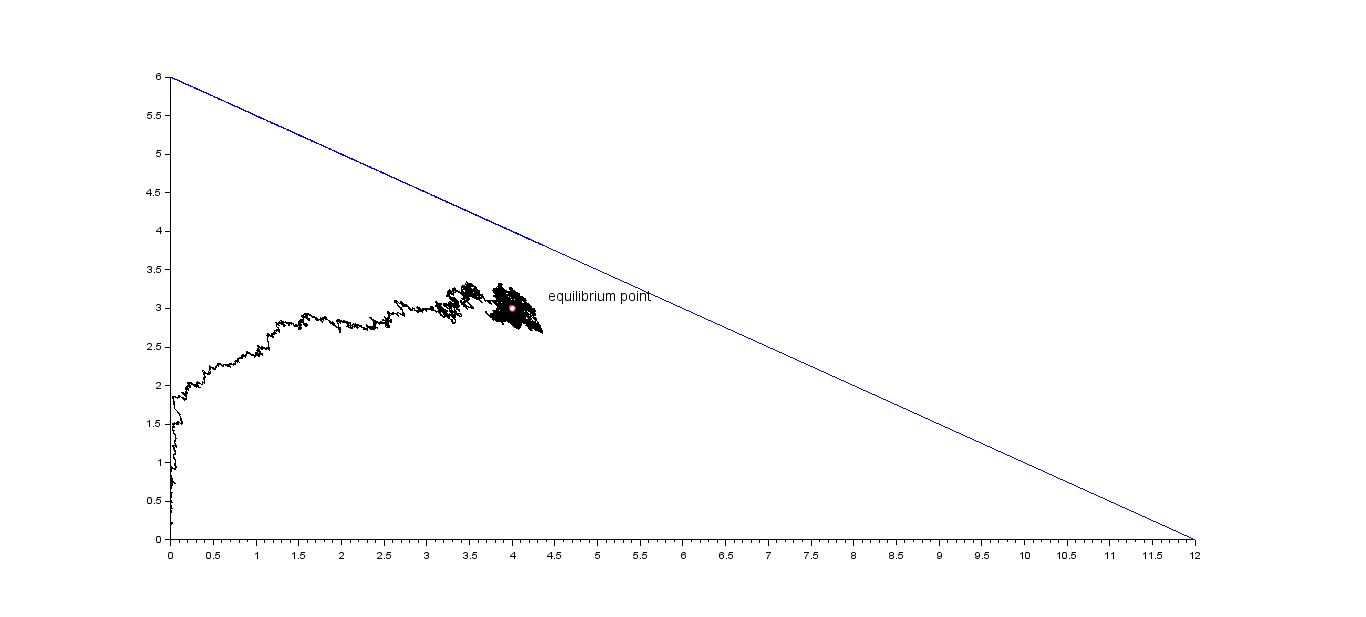} 
\caption{The under-load regime}
\label{fig3}
\end{figure}
Let $(X_1^N(t))$ and $(X_2^N(t))$ the processes given respectively by equations \eqref{eq4} and \eqref{eq5}.Recall that $(X_1^N(t)+X_2^N(t))$ is the process describing the total number of files that are present in the system at time $t$. Let $(Z^N(t))$  be the process given by 
\begin{equation}\label{24}
Z^N(t)=\dfrac{X_1^N(t)+2X_2^N(t)-N\rho}{\sqrt{N}}
\end{equation}

The  Q-matrix $Q^{N} = (q^{N}(.,.))$ of the Markov process $ (X_1^N(t) ,X_2^N(t) ,Z^N(t)) $ is defined by;

For $(x_1,x_2)\in\mathcal{D}^N$ and $z=\dfrac{x_1+2x_2-N\rho}{\sqrt{N}}$
\begin{equation}\label{25} (x_1,x_2,z) \longrightarrow (x_1,x_2,z)+
      \left\{  \begin{array}{ll}
       (0,1,\frac{2}{\sqrt{N}}) \ \  \xi N  \mathds{1}_{\lbrace z < \frac{F_{N} -1}{\sqrt{N}}-N\rho   \rbrace} \\
       (1,-1,-\frac{1}{\sqrt{N}})\ \  2 \mu  x_2\\
       (-1,1,\frac{1}{\sqrt{N}}) \ \  \lambda N   \mathds{1}_{\lbrace x_1> 0,   z < \frac{F_{N}}{\sqrt{N}}-N\rho   \rbrace}\\
       (-1,0,-\frac{1}{\sqrt{N}}) \ \ \mu x_1
      \end{array}
              \right.
  \end{equation}
and the generator of $ (X_1^N(t) ,X_2^N(t) ,Z^N(t)) $ is given by,
  \begin{align}\nonumber 
  A^N f(x_1, x_2, z) & = \xi N   \mathds{1}_{\lbrace z < \frac{F_N - N \rho - 1}{\sqrt{N}} \rbrace}  [f(x_1 , x_2 +1, z+ \frac{2}{\sqrt{N}}) - f(x_1, x_2, z)] \\\nonumber
&+ \lambda N \mathds{1}_{\lbrace x_1 > 0 \ , \ z < \frac{F_N - N \rho}{\sqrt{N}} \rbrace} [f(x_1 - 1, x_2 + 1, z+ \frac{1}{\sqrt{N}})- f(x_1, x_2, z)] \\\nonumber & + \mu x_1 [f(x_1 - 1, x_2, z- \frac{1}{\sqrt{N}}) - f(x_1, x_2, z)] \\\nonumber
&+ 2 \mu x_2  [f(x_1 + 1, x_2 - 1, z- \frac{1} {\sqrt{N}}) - f(x_1, x_2, z)] \end{align}
For any function $f$ depending only on the third variable $z$ i.e $$f(x_1,x_2,z)=g(z)\quad \forall\;(x_1,x_2)\in\mathbb{N}^2 \quad \textit{with}\quad x_1>0$$ for some twice differentiable function $g$ on $\mathbb{R}$ one gets
 \begin{align}\nonumber 
  A^N g( z) & = \xi N   \mathds{1}_{\lbrace z < \frac{F_N - N \rho - 1}{\sqrt{N}} \rbrace}  [g( z+ \frac{2}{\sqrt{N}}) - g( z)] \\\nonumber
&+ \lambda N \mathds{1}_{\lbrace  z < \frac{F_N - N \rho}{\sqrt{N}} \rbrace} [g( z+ \frac{1}{\sqrt{N}})- g( z)] \\\nonumber & + \mu(\sqrt{N}z+N\rho)[g(z- \frac{1}{\sqrt{N}}) - g( z)]  \end{align}
Remark that condition \eqref{23} implies that  terms $\frac{F_N - N \rho - 1}{\sqrt{N}}$ and $\frac{F_N - N \rho}{\sqrt{N}}$  converge to $+\infty$. Thus the generator converges to 
$$-\mu zg'(z)+(\lambda+3\xi)g''(z)\qquad z\in\mathbb{R}$$ 
when $N\rightarrow +\infty$, which is the generator of an Ornstein-Uhlenbeck process with variance converges to $\frac{\lambda+3\xi}{\mu}$. By results given in \cite{EK} one can see that for some positive constant $\alpha$ the process $(X_1^N(t)+2X_2^N(t))$ lives in $[N\rho-\alpha N, N\rho+\alpha N]\subset [0,N\bar{\beta}]$ and  the probability of saturation of the system is therefore small. In the under-loaded regime one can suppose that the capacity of the system is infinite i.e $F_N=+\infty$. In this case the complete study of the process $(X_1^N(t),X_2^N(t))$  is made in the article  \cite{EL}.

\section{The overloaded regime}

Throughout this section, we assume that the condition
\begin{equation}\label{24}
\rho>\bar{\beta}
\end{equation}
holds. 
\begin{figure}[hbtp]
\centering
\includegraphics[scale=0.5]{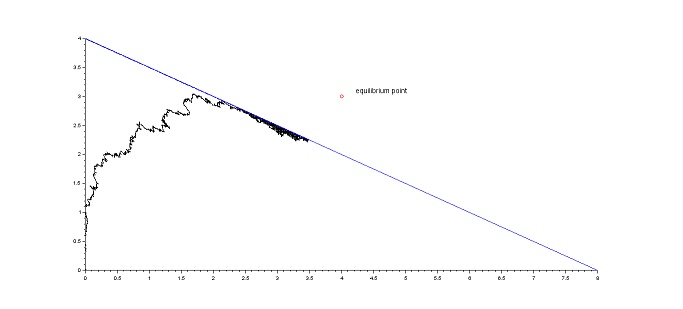}
\caption{The overload regime}
\label{fig4}
\end{figure}

the Q-matrix $Q^{N} = (q^{N}(.,.))$ and  the generator of the Markov  process $ ({X_1}^{N}(t/N), {X_2}^{N}(t/N) ,m^N(t/N)) $ are given by,
      \begin{equation}\nonumber 
      \left\{  \begin{array}{ll}
       q^{N}((x_1,x_2, m),(x_1 - 1, x_2, m+1)=  \frac{1}{N} \mu   x_1\\
        q^{N}((x_1,x_2, m),(x_1 +1, x_2-1, m+ 1)=  \frac{2}{N}  \mu x_2  \\
        q^{N}((x_1,x_2, m),(x_1 -1, x_2 + 1, m- 1)= \lambda \mathds{1}_{\lbrace  x_1 > 0 \ , \ m\geq 1 \rbrace} \\
        q^{N}((x_1,x_2, m),(x_1 , x_2 +1, m- 2)=  \xi  \mathds{1}_{\lbrace m\geq 2  \rbrace} 
      \end{array}
              \right.
  \end{equation}
 
  \begin{align}\nonumber 
  A_N f(x_1, x_2, m) & =  \frac{1}{N} \mu   x_1   [f(x_1 - 1, x_2, m+1) - f(x_1, x_2, m)] \\\nonumber
&+  \frac{2}{N}  \mu x_2 [f(x_1 +1, x_2-1, m+ 1)- f(x_1, x_2, m)] \\\nonumber & + \lambda \mathds{1}_{\lbrace  x_1 > 0 \ , \ F_N  - (2  x_2 +  x_1) \geq 1 \rbrace} [f(x_1 -1, x_2 + 1, m- 1) - f(x_1, x_2, m)] \\\nonumber
&+ \xi  \mathds{1}_{\lbrace F_N  - (2  x_2 +  x_1) \geq 2  \rbrace}   [f(x_1 , x_2 +1, m- 2) - f(x_1, x_2, m)]  
 \end{align}
For any function $f$ depending only on the third variable $m$ i.e $$f(x_1,x_2,m)=g(m)\quad \forall\;(x_1,x_2)\in\mathbb{N}^2 $$ for some function $g$ on $\mathbb{N}$ one gets
\begin{align}\nonumber 
  A_N g( m) & =    \mu\frac{F_N-m}{N}\left(g(m+1) - g(m)\right) \\\nonumber
&+ \lambda \mathds{1}_{\lbrace  x_1 > 0 \ ,  m \geq 1 \rbrace} [g( m- 1) - g(m)] \\\nonumber
&+ \xi  \mathds{1}_{\lbrace m \geq 2  \rbrace}   [g(m- 2) - g( m)]  
 \end{align}
 This generator converges to 
 \begin{align}\nonumber 
  Ag( m) & =    \mu\bar{\beta}\left(g(m+1) - g(m)\right) \\\nonumber
&+ \lambda \mathds{1}_{\lbrace x_1>0, m \geq 1 \rbrace} [g( m- 1) - g(m)] \\\nonumber
&+ \xi  \mathds{1}_{\lbrace m \geq 2  \rbrace}   [g(m- 2) - g( m)]  
 \end{align} 
 Thus, for any $x=(x_1,x_2)\in\mathbb{N}^*\times\mathbb{N}$,  this is the  generator of the Markov process $(m(t))$ with transitions
 \begin{equation}\label{qmatrix} m \longrightarrow m+
      \left\{  \begin{array}{ll}
      +1 \ \  \mu \bar{\beta} \\
      -1 \ \  \lambda \mathds{1}_{\{ m \geq 1\}}  \\
      -2 \ \  \xi \mathds{1}_{\{ m \geq 2\}} 
      \end{array}
              \right.
  \end{equation}
 \begin{proposition} Under the condition \eqref{24}, the process $(m(t))$ has a unique invariant distribution $\pi$, its generating function $g(u) = {\underset{n \geq 0}{\sum }} \pi(n) u^n$ is given by, for $u\in[-1,1]$
 \begin{equation}\label{25}
  g(u) = \frac{1}{-\mu \bar{\beta}+ (\lambda+\xi) u +\xi} \left[(\lambda u+\xi (1+u)) \pi(0) +\xi (1+u) u  \pi(1)\right]
 \end{equation}

Where $( \pi(0), \pi(1))$ are given by
\begin{equation}\label{pi0}
\pi(0) = \frac{(1+y_{*}) (\lambda+2 \xi - \mu \bar{\beta})}{(\lambda+2\xi)( 1+y_{*}) -2\mu\bar{\beta}y_{*}}
\end{equation}
\begin{equation}\label{pi1}
\pi(1)= \frac{-\mu \bar{\beta} + \lambda +2 \xi}{2 \xi} - \frac{\lambda + 2 \xi}{2 \xi }  \pi(0)
\end{equation}

with $$y_{*} = \frac{(\lambda+\xi)- \sqrt{(\lambda+\xi)^2 + 4 \xi \mu \bar{\beta}}}{2 \mu \bar{\beta}}$$
\end{proposition}\label{prop1}
\begin{proof}
Existence and uniqueness of the stationary distribution is a simple consequence of Foster's criterion. See Proposition 8.14 of \cite{robert2}. For $u \in [-1,1]$, define 
$$g(u) = {\underset{n \geq 0}{\sum }} \pi(n) u^n$$
The equilibrium equation 
\begin{align}\label{26}
&\sum_{m=0}^{+\infty}[\mu\bar{\beta}\left(f(m+1)-f(m)\right)+\lambda\mathds{1}_{\lbrace  x_1 > 0 \ ,  m \geq 1 \rbrace}\left(f(m-1)-f(m)\right)\\\nonumber
&\qquad+\xi\mathds{1}_{\lbrace   m \geq 2 \rbrace}\left(f(m-2)-f(m)\right)] \pi(m)=0
\end{align}

for $f(m)=u^m$, gives the following relation

\begin{equation*}
g(u) (\mu \bar{\beta} u^2 (u-1) + \lambda (u-u^2) +\xi (1-u^2)) = \lambda (u-u^2)  \pi(0) + \xi (1-u^2) ( \pi(0) + u  \pi(1))
\end{equation*}

Let $$P(u)\stackrel{def}{=} - \mu \bar{\beta} u^2 +(\lambda +\xi) u +\xi $$
then we have 
\begin{equation}\label{27}
P(u) g(u) = ((\lambda+\xi) u+ \xi )\pi (0) + \xi (1+u) u \pi(1))
\end{equation}
Note that $P(-1)=-( \mu \bar{\beta}+\lambda)<0$,\;$P(0)=\xi$ and  $P(1)=-\mu \bar{\beta}+\lambda+2\xi >0$ by Condition \eqref{24}. The function $P(u)$ has a unique root in$[-1,1]$ and it is necessarily $y_*$.

We have therefore  that $y_*$ is a root of the right-hand side of the Relation \eqref{27}, hence 

\begin{equation*}
\mu\bar{\beta}y_*^2\pi(0)+\xi y_*(1+y_*)\pi(1)=0
\end{equation*}
 
and the relation $g(1)=1$ gives the additional identity 
\begin{equation*}
\frac{\lambda+2\xi}{2 \xi} \pi(0)+\pi(1) =\frac{\lambda + 2 \xi - \mu \bar{\beta}}{2 \xi}
\end{equation*}  
The proposition is proved

\end{proof}
\subsection{Fluid limits}
Our aim in this section is to identify the limit  of the renormalized  processes  $(\bar{X}_1^N(t))$ and $(\bar{X}_2^N(t))$ given respectively by equations \eqref{eq11} and \eqref{12}. We assume that 

\begin{equation}\label{28}
\underset{N \rightarrow + \infty}{lim} (\overline{X}_1^N(0),\overline{X}_2^N(0)) =(x_1,x_2)
\end{equation}
and we successively study the cases where $(x_1,x_2)$ is chosen inside the set $\mathcal{S}$ and the case where $(x_1,x_2)$ lies on the boundary $(\partial \mathcal{S})_2$ .
\subsubsection{Starting from the interior of $ \mathcal{S}$}.

Let $T_1^N$ the hitting time   $$T_1^N=\inf\{t>0\;|\; m^N(t)\in\{0,1\}\}$$

Note that before time $T_1^N$ the Markov process $(X_1^N(t),X_2^N(t))$ coincides with  the Markov process  describing the storage process with infinite capacity ($F_N=+\infty$). 

The  Proposition~\ref{prop3}   proves the convergence in distribution of the hitting 
time $T_1^N$. The proof of this result is inspired by the study of $M/M/N/N$ queue ( see\cite{robert2} and \cite{CFR}). 
Let $\phi_c^N$ be the function on $\mathbb{R}^+$ defined by

$$\phi_c(t) = ce^{\mu t} (\rho + \frac{c \xi }{2 \mu} e^{\mu t})$$
for $c \in \mathbb{R}^*$,\;$N\in\mathbb{N}^*$. 
  
\begin{lemma}\label{lem2}Let $v=(1,2)$.  The function 
\begin{align*}
g_c : &(t,w) \in \mathbb{R}^+ \times \mathbb{N}^*\times\mathbb{N} \rightarrow (1+c e^{\mu t} )^ {v\cdot w} e^{- N\phi_c(t)}\\
&\text{where}\qquad v\cdot w=w_1+2w_2 
\end{align*}
is space-time harmonic with respect to the Q-matrix $Q$ given in \eqref{eq1} with $F_N=+\infty$. In other words 
$$\frac{\partial g_c}{\partial t} (t,w) + Q (g_c) (t,w) =0,\; \text{for all} \ t \in \mathbb{R}^+ \; \text{and  for all} \; w \in \mathbb{N}^*\times\mathbb{N}$$

\end{lemma}
\begin{proof}
For  $t\in\mathbb{R}_+$ and $w\in\mathbb{N}^*\times\mathbb{N}$
\begin{equation}\nonumber
\begin{split}
\frac{\partial g_c}{\partial t} (t,x) & =e^{- \phi_c (t)}ce^{\mu t} \biggl[  v\cdot w   \mu(1+c e^{\mu t})^ {v\cdot w-1}\\
 & \hspace{2cm}-  \biggl( \lambda N + \xi N(2+ c e^{\mu t})\biggr) (1+c e^{\mu t})^{v\cdot w}\biggr]
\end{split}
\end{equation} 
on other hand
$$Q(g_c)(t,w)=Q(g_c(t,.))(w)$$ is given by

\begin{equation}\nonumber 
\begin{split}
Q(g_c)(t,w)& =\lambda N \biggl[(1+c e^{\mu t})^{v\cdot w+1}  e^{- \phi_c (t)} - (1+c e^{\mu t})^{v\cdot w}  e^{- \phi_c (t)}\biggr]  \\
&\qquad + \mu  {v\cdot w}  \biggl[(1+c e^{\mu t})^{v\cdot w-1}  e^{- \phi_c (t)} - (1+c e^{\mu t})^{v\cdot w}  e^{- \phi_c (t)}]\\
& \quad\quad\quad+ \xi N \biggl[(1+c e^{\mu t})^{v\cdot w+2}  e^{- \phi_c (t)} - (1+c e^{\mu t})^{v\cdot w}e^{- \phi_c (t)}\biggr]\\
&\\
& = e^{- \phi_c (t)} \biggl[(\lambda + \xi) N  \biggl((1+c e^{\mu t})^{v\cdot w+1}  - (1+c e^{\mu t})^{v\cdot w} \biggr)\\
&\qquad +\mu{v\cdot w}  \biggl((1+c e^{\mu t})^{v\cdot w-1} - (1+c e^{\mu t})^{v\cdot w}\biggr )\\
& \qquad\qquad+  \xi N  \biggl((1+c e^{\mu t})^{v\cdot w+2} -(1+c e^{\mu t})^{v\cdot w+1}\biggr )\biggr]\\
&\\
& = e^{- \phi_c (t)} \biggl[- {v\cdot w}c \mu e^{\mu t} (1+c e^{\mu t})^{v\cdot w-1}  \\
&\qquad +c e^{\mu t} \biggl(\lambda N + \xi N (2 + c e^{\mu t})\biggr) (1+c e^{\mu t})^{v\cdot w}  \biggr]\\
&\qquad\qquad = - \frac{\partial g_c}{\partial t} (t,w)
 \end{split} 
\end{equation}
\end{proof}
\newpage
\begin{proposition}\label{prop2}  
 \begin{enumerate}
 \item For $c \in \mathbb{R}^*$ and $N\in\mathbb{N}^*$ the process
\begin{equation} \label{eq29}  \left( g_c(t,X^N(t))\right) 
\end{equation}
is a martingale.
\item For $N\in\mathbb{N}^*$ the following processes are martingales.
  \begin{equation} \label{eq30}\left(e^{\mu t} (v\cdot X^N(t) - N\rho )\right)
   \end{equation} 
   \begin{equation}\label{eq31} \left( e^{2 \mu t} \left((v\cdot X^N(t) - N\rho )^2 - v\cdot X^N(t) - N\frac{\xi}{\mu}\right)\right)
   \end{equation}
   \end{enumerate}
\end{proposition}
\begin{proof}
\begin{enumerate}
\item By Lemma~\ref{lem2} the function  $(t,w)\rightarrow  g_c(t,w)$ is space-time harmonic for the Q-matrix $Q$ given in (\eqref{eq1}) with $F_N=+\infty$. Since  $t \rightarrow \frac{\partial g_c}{\partial t}$ is continuous, then the process $(g_c(t,X^N(t))$ is  a local martingale (See Corollary B.5 in   \cite{robert2}).  Furthermore, for $t \in \mathbb{R}^+$,  
 \begin{eqnarray} \nonumber v\cdot X^N(t)  \leq (2 X_2^N(0) + X_1^N(0))+ 2 \mathcal{N}_{\xi_N}(]0,t]) + \mathcal{N}_{\lambda_N}(]0,t]) 
 \end{eqnarray}
  one gets for $t\geq 0$,  \[\mathbb{E} ( \sup_{0 \leq s \leq t }|g_c(t,X^N(t)|) < + \infty\]
 Thus the  process $(g_c(t,X^N(t))$ is a martingale (see proposition A.7 in \cite{robert2}).
 \item Let $\Psi$ be the function  on $\mathbb{R}^+\times \mathbb{N}$ defined  by

 \[\Psi(x,z)=(1+x)^ze^{-N\rho x}e^{-\frac{N\xi}{2\mu}x^2}\]
 Note that  \[\Psi(ce^{\mu t},v\cdot X^N(t))=g_c(t,X^N(t))\] and therefore $(\Psi(ce^{\mu t},v\cdot X^N(t))$ is a martingale. On other hand it is well known that 
 \[e^{-N\rho x}(1+x)^z=\sum_{n\geq 0}C_n^{N\rho}(z)\frac{x^n}{n!}
 \] where $C_n^{N\rho}(z)$ is the $n$th Poisson-Charlier polynomial (see \cite{chi}) . Hence the  expansion of $\Psi(x,z)$ is given by
 
 \begin{equation}\label{32}
 \Psi(x,z)=\sum_{n\geq 0}\left(\sum_{k=0}^nC_{n-k}^{N\rho}(z)b_k\right)\frac{x^n}{n!}
 \end{equation}
  where $b_{2k+1}=0$ and $b_{2k}=\left(-\frac{N\xi}{2\mu}\right)^k$ 
 \end{enumerate}
 Replacing  in  \eqref{32} $x$ and $z$ by $ce^{\mu t}$ and  $v\cdot X^N(t)$ respectively one gets  that for any $n\in\mathbb{N}^*$, 
 \[\left( e^{n\mu t}(\sum_{k=0}^nC_{n-k}^{N\rho}(v\cdot X^N(t))b_k)\right)\]
 is a martingale. In particular for $n=1$ and $n=2$ one gets that the processes 
 \[\left(e^{\mu t}(v\cdot X^N(t)-N\rho)\right)\]
 and 
 \[\left(e^{2\mu t}(v\cdot X^N(t)-N\rho)^2-v\cdot X^N(t)-\frac{N\xi}{2\mu})\right)\]
 are martingales.
\end{proof}
\begin{proposition} \label{prop3}

if Conditions \eqref{24} and \eqref{28} hold with $x_1+2x_2<\bar{\beta}$  then the hitting time $T_{1}^N$ converges  in distribution to $T_0$ where
\begin{equation} T_0 = \frac{1}{\mu} log(\frac{\lambda+ 2 \xi - \mu (x_1+2 x_2)}{\lambda+ 2 \xi - \mu \bar{\beta}})
\end{equation}

\end{proposition}
\begin{proof}
We assume that Conditions \eqref{24} and \eqref{28} hold with $x_1+2x_2<\bar{\beta}$.
Doob's optional stopping Theorem applied to the martingale given in \eqref{eq30} and to $T_1^N$  show that the  process 
 \begin{equation} \nonumber
 \left(e^{\mu t \wedge T_{1}^N}\left[ v\cdot X^{N}(t \wedge T_{1}^N)- N\rho \right]\right)
 \end{equation}
is a martingale. Thus  the following equality holds
\begin{equation*}  \mathbb{E}  \left(e^{\mu t \wedge T_{1}^N}\left[N\rho  -v\cdot X^{N}(t \wedge T_{1}^N) \right]\right)=  N\rho-v\cdot X^{N}(0) 
\end{equation*}
 Since $v\cdot X^{N} (t \wedge T_{1}^N) \leq F_N-1$, one gets   that, 
\begin{equation} \nonumber \mathbb{E} (e^{\mu  t \wedge T_{1}^N} ) \leq \frac{(\lambda + 2\xi) N - \mu v\cdot X^{N}(0)}{(\lambda + 2\xi) N - \mu F_N+\mu}
 \end{equation}
By letting $t$ go to infinity, monotone convergence Theorem shows that 
 \begin{equation}  \nonumber \mathbb{E} (e^{\mu  T_{1}^N} ) \leq \frac{\lambda + 2\xi - \mu v\cdot\bar{X}^N(0)}{\lambda + 2\xi - \mu \frac{F_N}{N}+\frac{\mu}{N}} 
 \end{equation}
And that implies uniform integrability of the martingale 
\begin{equation} \nonumber \mathbb{E} (e^{\mu  t \wedge T_{1}^N} (v\cdot X^{N} (t \wedge T_{1}^N)- \rho N))
\end{equation}
One gets therefore the following identity 
\begin{equation} \label{eq32}\mathbb{E} (e^{\mu T_1^N}) = \frac{ \lambda +2 \xi -\mu v\cdot\bar{X}^N(0)}{ \lambda + 2 \xi - \mu \frac{\bar{F}_N}{N}+\frac{\mu}{N} } 
\end{equation}
Doob's optional stopping theorem applied again to the martingale given by \eqref{eq31} and to the stopping time $T_1^N$ shows that the process  
\begin{equation}  \nonumber (e^{{2 \mu t \wedge  T_{1}^N}} (v\cdot X^{N}(t \wedge  T_{1}^N)- \rho N)^{2} - v\cdot X^{N}(t\wedge  T_{1}^N) - \frac{\xi N}{\mu}))
\end{equation} is a martingale. Since $v\cdot X^{N} (t \wedge T_{1}^N) \leq F_N-1$,$N\bar{\beta}<N\rho$  and $N\bar{\beta}= F_N$  one could then use the same arguments used above to get the following identity

\begin{equation} \label{eq33}\mathbb{E} (e^{{2 \mu T_{1}^N}}) = \frac{N (v\cdot\bar{X}^N(0)- \rho )^{2} - v\cdot\bar{X}^N(0) - \frac{\xi }{\mu}}{N(\frac{F_N-1}{N}- \rho)^{2} - \frac{F_N-1}{N} - \frac{\xi }{\mu}}
\end{equation}
One then deduces that $var(e^{\mu T_{1}^N}) = O(1/N)$ and the Tchebychev inequality implies that, for $\epsilon > 0$, 
 $$\mathbb{P} (| e^{\mu  T_{1}^N}  - \mathbb{E} (e^{\mu  T_{1}^N} ) | > \epsilon) \leq  \frac{var(e^{\mu  T_{1}^N})}{\epsilon^{2}},$$
 Hence, using the identity given by \eqref{eq32}, the sequence $(T_1^N)$ converges in probability to $T_0$. 
\end{proof}
\begin{theorem}\label{th2} If Conditions \eqref{24} and \eqref{28} hold with $x_1+2x_2<\bar{\beta}$ and $x_2>\frac{\lambda+\xi}{2\mu}$
Then for the convergence in distribution, 
\begin{equation} \nonumber \underset{N \rightarrow + \infty}{lim} (\overline{X}_1^N(t),\overline{X}_2^N(t))_{0\leq t \leq T_0} =(\bar{x}_1(t),\bar{x}_2(t))_{0\leq t \leq T_0}  
\end{equation} with $(\bar{x}_1(t), \bar{x}_2(t)) $ are given in \eqref{22}.

Note that at time $T_0$ , the fluid limit $(\bar{x}_1(t), \bar{x}_2(t)) $ hits the boundary $(\partial\mathcal{S})_2$ i.e $\bar{x}_1(T_0) + 2 \bar{x}_2(T_0) = \bar{\beta}$.

\end{theorem}

\begin{proof}
We assume that  Conditions \eqref{24} and \eqref{28} hold with $x_1+2x_2<\bar{\beta}$ and $x_2>\frac{\lambda+\xi}{2\mu}$. 
By Theorem~\ref{th1} the sequence \[\left(\overline{X}^{N}(t),Y^{N}(t), \nu^N(t)\right)\] is relatively compact in $\mathcal{D}(\mathbb{R}_+,\mathbb{R}^3)$ and the limit $\left( x(.) ,y(.), \nu(.)\right)$ of any convergent subsequence satisfies for all $t\geq 0$ : 
\begin{align}\label{34}x_1(t)  &=  x_1 - \lambda t-\mu \int_{0}^t x_1(s)ds  + 2\mu \int_{0}^t x_2(s) ds \\\nonumber
   &+ \lambda \int_{[0,t] \times \mathbb{N}} \mathds{1}_{\lbrace x_1(s) > 0 \rbrace} \mathds{1}_{\lbrace  0 \rbrace} (u) \nu (ds \times du)   +\lambda y_1(t)
 \end{align}
\begin{align} \label{35}  x_2(t)  &=  x_2 + (\lambda+\xi) t -2  \mu \int_{0}^t x_2(s)  ds- \xi  \nu ([0,t] \times \{0,1\})  \\\nonumber
   &  - \lambda \int_{[0,t] \times  \mathbb{N}}  \mathds{1}_{\lbrace x_1(s) > 0 \rbrace} \mathds{1}_{\lbrace  0 \rbrace} (u)  \nu (ds \times du) - \lambda y_1(t)  
 \end{align}
 The condition  $x_2>\frac{\lambda+\xi}{2\mu}$ implies that the function $y_1(t)=0$ for all $t\geq 0$  (see Theorem 2 in \cite{EL}). Thus, it is sufficient to show that for all $t\leq T_0$
\[\nu ([0,t] \times \{0,1\})=0\] 
 
Let us first recall that,\[ \nu^N((0,t)\times \{0,1\})=\int_0^t\mathds{1}_{\{m^N(u)\in\{0,1\}\}}du\]
and that the increasing  sequence of hitting times $(T_1^N)$ converges in probability to $T_0$. For  any  $t\leq T_0$ and for any $\epsilon >0$
\begin{align*}
\mathbb{P}\{\sup_{s\leq t}\nu^N((0,s)\times \{0,1\})\geq\epsilon\}&\leq
\mathbb{P}\{\sup_{s\leq t\wedge T_1^N }\nu^N((0,s)\times \{0,1\})\geq\epsilon\}\\
\qquad &+ \mathbb{P}\{\sup_{T_1^N\leq s\leq t  }\nu^N((0,s)\times \{0,1\})\geq\epsilon\}
\end{align*}

The first term of the right-hand side of the above Inequality  is equal to zero. Since for $T_1^N\leq s\leq t $ 
\begin{align*}
\nu^N((0,s)\times \{0,1\})&=\int_0^{ T_1^N}\mathds{1}_{\{m^N(u)\in\{0,1\}\}}du+\int_{T_1^N}^s\mathds{1}_{\{m^N(u)\in\{0,1\}\}}du\\
&\leq T_0-T_1^N
\end{align*} 
 \[\mathbb{P}\{\sup_{s\leq t}\nu^N((0,s)\times \{0,1\})\geq\epsilon\}\leq
\mathbb{P}\{|T^N-T_0|\geq\epsilon\} \]
Thus, \[\lim_{N\rightarrow +\infty}\mathbb{P}\{\sup_{s\leq t}\nu^N((0,s)\times \{0,1\})\geq\epsilon\}=0\]
 \end{proof}

\subsubsection{Starting from the boundary of $(\partial\mathcal{S})_2$ of the set $\mathcal{S}$}

\begin{theorem}\label{3}  If Conditions \eqref{24} and \eqref{28} hold with $x_1+2x_2=\bar{\beta}$.  Then for the convergence in distribution,

 \[   \underset{N \rightarrow + \infty}{lim}  (\bar{X}_1^N(t),\bar{X}_2^N(t))_{t \geq 0}   = (x_1(t), x_2(t))_{t \geq 0}  
 \]
where $(x_1(t), x_2(t))_{t \geq 0}  $ is the solution of the ordinary differential equation,
\begin{align} x_1(t)  &=  x_1 - \lambda (1-\pi(0))t + \mu \int_{0}^t (2 x_2(u) - x_1(u)) du \\\nonumber
  x_2(t)  &=  x_2 +   \left(\frac{\mu\bar{\beta}}{2}+\frac{\lambda}{2}(1-\pi(0)\right) t -2  \mu \int_{0}^t x_2(u)  du 
 \end{align}
 where $\pi(0)$ is defined by Equation \eqref{pi0}.

\end{theorem} 

\begin{proof}
Our goal is to identify   the measure $\nu$ in Equations \eqref{34} and \eqref{35}. The  Q-matrix of the Markov process $(X^N(.), m^N(.))$ is given by,
 \begin{equation}\nonumber (x^N,m^N) \longrightarrow (x^N,m^N)+
      \left\{  \begin{array}{ll}
       (x^N+e_2,m^N-2) \ \  \xi N \mathds{1}_{\{m^N \geq 2\}} \\
       (x^N+e_1-e_2,m^N+1)\ \  2 \mu   x_2^N\\
       (x^N+e_2-e_1,m^N-1) \ \  \lambda N  \mathds{1}_{\{ x_1^N > 0\}} \mathds{1}_{\{ m^N \geq 1\}} \\
       (x^N-\frac{e_1}{N},m^N+1) \ \ \mu x_1^N
      \end{array}
              \right.
  \end{equation}
  Thus,the process 
  \[\left(f(X^N(t), m^N(t) ) - f(X^N(0), m^N(0) )-\int_0^t (Qf)(X^N(s), m^N(s) )ds\right)
\]
is a  martingale for all bounded function $f$ on $\mathbb{R}^+\times \bar{\mathbb{N}}$. In particular 
 the process 
\begin{equation}\label{36}
\begin{split}
 \mathcal{M}^N(t)\stackrel{def}{=}&g(m^N(t)) - g(m^N(0))\\
 & - \int_{0}^t [g(m^N(s)-2)- g(m^N(s))]\xi N \mathds{1}_{\{m^N(s) \geq 2\}}ds \\
&\quad -\int_{0}^t  [g(m^N(s)+1) - g(m^N(s))]  \mu(2X_2^N+X_1^N(s))ds \\
&\qquad -\int_{0}^t [g(m^N(s)-1) - g(m^N(s))]  \lambda N  \mathds{1}_{\{X_1^N(s) > 0\}} \mathds{1}_{\{ m^N(s) \geq 1\}}ds
\end{split}
\end{equation}
is a martingale for all bounded function $g$ on $\bar{\mathbb{N}}$. It follows from  Doob's inequality  that the process $ (\frac{\mathcal{M}^N(t)}{N})$ converges in distribution to 0.

 Since $2\bar{X}_2^N(t)+\bar{X}_1^N(t)=\frac{F_N}{N}-\frac{m^N(t)}{N}$, Equation \eqref{36} can be rewritten as 
\begin{equation}\label{37}
\begin{split}
 \frac{\mathcal{M}^N(t)}{N}=&\frac{g(m^N(t)) - g(m^N(0))}{N}\\
 & - \int_{0}^t\biggl\{ [g(m^N(s)-2)- g(m^N(s))]\xi \mathds{1}_{\{m^N(s) \geq 2\}} \\
&\quad + [g(m^N(s)+1) - g(m^N(s))]  \mu(\frac{F_N}{N}-\frac{m^N(t)}{N}) \\
&\qquad + [g(m^N(s)-1) - g(m^N(s))]  \lambda  \mathds{1}_{\{X_1^N(s) > 0\}} \mathds{1}_{\{ m^N(s) \geq 1\}}\biggr\} ds
\end{split}
\end{equation}

In terms of measure $\nu^N(.)$, we may rewrite the last term on the right-hand side of \eqref{37} as follows :
  
\begin{equation}\label{38}
\begin{split}
 &  \int_{0}^t\biggl\{ (g(y-2)- g(y))\xi \mathds{1}_{\{y \geq 2\}}\\
&\quad +  (g(y+1) - g(y))  \mu(\frac{F_N}{N}-\frac{y}{N})\\
&\qquad + [g(y-1) - g(y)]  \lambda  \mathds{1}_{\{\bar{X}_1^N(s) > 0\}} \mathds{1}_{\{ y \geq 1\}}\biggr\}\nu^N(ds\times dy)
\end{split}
\end{equation}
which, converges  also to 0 since  $\frac{1}{N}\left(g(m^N(t)) - g(m^N(0))\right)$ converges to 0
 as $N\rightarrow+\infty$.   Furthermore by continuous mapping theorem one gets that    
\begin{align}\nonumber   \int_{[0,t]\times \mathbb{N}}  \biggl\{ [g(y-2)- g(y)] \xi \mathds{1}_{\{y \geq 2\}} +  [g(y+1) - g(y)]  \mu \bar{\beta}  \\\nonumber
+ [g(y-1) - g(y)] \lambda  \mathds{1}_{\{ x_1(s) > 0\}} \mathds{1}_{\{ y \geq 1\}} \biggr\} \nu(ds \times dy) =0
 \end{align}
 for all $t\geq 0$.
 
 Thus for almost all $t \geq 0$, 
 \begin{align}\nonumber   \sum_{y\in \mathbb{N}}  \biggl\{ [g(y-2)- g(y)] \xi \mathds{1}_{\{y \geq 2\}} +  [g(y+1) - g(y)]  \mu \bar{\beta}  \\\nonumber
+ [g(y-1) - g(y)] \lambda  \mathds{1}_{\{x_1(t) > 0\}} \mathds{1}_{\{ y \geq 1\}} \biggr\}\nu_t (y) =0
 \end{align}
 Hence, for all $t\geq 0$ such that $x_1(t)>0$ the measure $\nu_t(.)=\pi$ where the measure $\pi$ is invariant for the Markov process $(m(t))$ with Q-matrix given by \eqref{qmatrix}.
 The Theorem is proved.
\end{proof}

\section{The critical regime}
\begin{figure}[hbtp]
\centering
\includegraphics[scale=0.2]{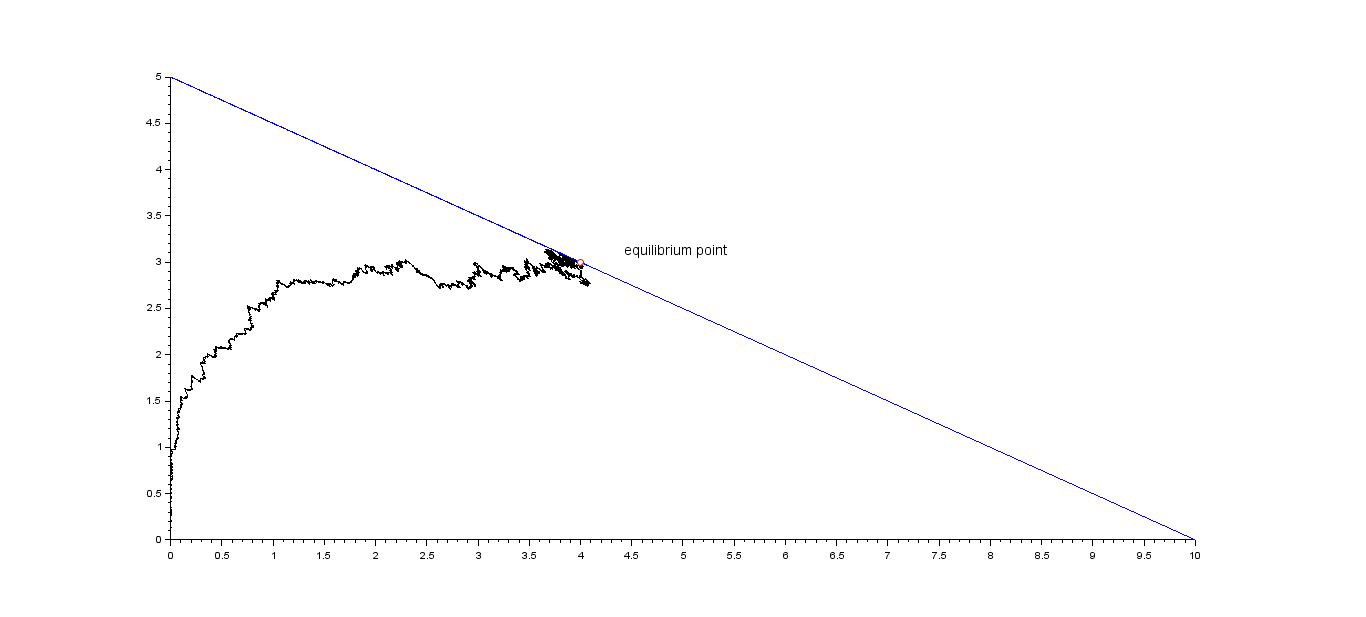}
\caption{The critical regime}
\label{fig2}
\end{figure}
Throughout this section, we assume that the condition
\begin{equation}\label{39}
\rho=\bar{\beta}
\end{equation}
holds.
Let $\lbrace Z_1^N(t),Z_2^N(t),Z^N(t)\rbrace$ be the Markov process defined by 
\[Z_1^N(t)=\sqrt{N}(\rho_1-\bar{X}_1^N(t)),\quad Z_2^N(t)=\sqrt{N}(\rho_2-\bar{X}_2^N(t))\]
and \[Z^N(t)=\sqrt{N}(\rho-\bar{X}_1^N(t)-2\bar{X}_2^N(t))=Z_1^N(t)+2Z_1^N(t)\]
 where \[\rho_1=\frac{\xi}{\mu},\quad \rho_2=\frac{\lambda+\xi}{2\mu}\]
 In the following proposition we prove that the sequence of processes \[\lbrace Z_1^N(t),Z_2^N(t),Z^N(t)\rbrace\] converges in distribution to a reflected three dimensional Ornstein-Uhlenbeck process.
 
The Q-matrix $Q^{N} = (q^{N}(.,.))$ and  the generator of the Markov  process $\lbrace Z_1^N(t),Z_2^N(t),Z^N(t)\rbrace$

are given by,

     \begin{equation}\nonumber 
      \left\{  \begin{array}{ll}
       q^{N}((z_1,z_2, z),(z_1+\frac{1}{\sqrt{N}}, z_2, z+\frac{1}{\sqrt{N}})=\mu N(\rho_1-\frac{z_1}{\sqrt{N}})\\
        q^{N}((z_1,z_2, z),(z_1-\frac{1}{\sqrt{N}}, z_2+\frac{1}{\sqrt{N}}, z+\frac{1}{\sqrt{N}})=  2\mu N(\rho_2-\frac{z_2}{\sqrt{N}} ) \\
        q^{N}((z_1,z_2, z),(z_1+\frac{1}{\sqrt{N}}, z_2-\frac{1}{\sqrt{N}}, z-\frac{1}{\sqrt{N}})= \lambda N \mathds{1}_{\lbrace  z_1 <\sqrt{N}\rho_1,\;z\geq \frac{1}{\sqrt{N}}+\sqrt{N}(\bar{\beta}-\frac{F_N}{N})   \rbrace} \\
       q^{N}((z_1,z_2, z),(z_1, z_2-\frac{1}{\sqrt{N}}, z-2\frac{1}{\sqrt{N}})= \xi N \mathds{1}_{\lbrace  z\geq \frac{2}{\sqrt{N}}+\sqrt{N}(\bar{\beta}-\frac{F_N}{N})   \rbrace}  
      \end{array}
              \right.
  \end{equation}

  \begin{align}\nonumber 
  &A_N f(z_1, z_2, z)  =  \mu N(\rho_1-\frac{z_1}{\sqrt{N}})  [f(z_1+\frac{1}{\sqrt{N}}, z_2, z+\frac{1}{\sqrt{N}}) - f(z_1,z_2, z)] \\\nonumber
&+ 2\mu N(\rho_2-\frac{z_2}{\sqrt{N}} ) [f(z_1-\frac{1}{\sqrt{N}}, z_2+\frac{1}{\sqrt{N}}, z+\frac{1}{\sqrt{N}})- f(z_1,z_2, z)] \\\nonumber
 & + \lambda N \mathds{1}_{\lbrace  z_1 <\sqrt{N}\rho_1,\;z\geq \frac{1}{\sqrt{N}}+\sqrt{N}(\bar{\beta}-\frac{F_N}{N})   \rbrace} [f(z_1+\frac{1}{\sqrt{N}}, z_2-\frac{1}{\sqrt{N}}, z-\frac{1}{\sqrt{N}}) - f(z_1,z_2, z)] \\\nonumber
&+ \xi N \mathds{1}_{\lbrace  z\geq \frac{2}{\sqrt{N}}+\sqrt{N}(\bar{\beta}-\frac{F_N}{N})   \rbrace}    [f(z_1, z_2-\frac{1}{\sqrt{N}}, z-2\frac{1}{\sqrt{N}}) - f(z_1,z_2, z)]  
 \end{align}
\begin{proposition}\label{prop4}
 If $f$ is twice differentiable on $\mathbb{R}^3$ and such that $\nabla f(z_1,z_2,0)=0$ then the generator  converges to 
\begin{equation}\label{•}
\begin{split}
Af(z_1,z_2,z)&=\mu (2z_2-z_1)\frac{\partial f}{\partial x_1}(z_1,z_2,z)-2\mu z_1\frac{\partial f}{\partial x_2}(z_1,z_2,z)\\
&-\mu z\frac{\partial f}{\partial x_3}(z_1,z_2,z)+(\lambda+\xi)\frac{\partial^2 f}{\partial x_1^2}(z_1,z_2,z)+(\lambda+\xi)\frac{\partial^2 f}{\partial x_2^2}(z_1,z_2,z)\\
&(\lambda+\frac{5}{2}\xi)\frac{\partial^2 f}{\partial x_3^2}(z_1,z_2,z)-(2\lambda+\xi)\frac{\partial^2 f}{\partial x_1\partial x_2}(z_1,z_2,z)-2\lambda\frac{\partial^2 f}{\partial x_1\partial x_3}(z_1,z_2,z)\\
&+(2\lambda+3\xi)\frac{\partial^2 f}{\partial x_2\partial x_3}(z_1,z_2,z)
\end{split}
\end{equation}
for $z>0$ and to 
\begin{equation}\label{•}
\begin{split}
&(\lambda+\xi)\frac{\partial^2 f}{\partial x_1^2}(z_1,z_2,0)+(\lambda+\xi)\frac{\partial^2 f}{\partial x_2^2}(z_1,z_2,0)\\
&(\lambda+\frac{5}{2}\xi)\frac{\partial^2 f}{\partial x_3^2}(z_1,z_2,0)-(2\lambda+\xi)\frac{\partial^2 f}{\partial x_1\partial x_2}(z_1,z_2,0)-2\lambda\frac{\partial^2 f}{\partial x_1\partial x_3}(z_1,z_2,0)\\
&+(2\lambda+3\xi)\frac{\partial^2 f}{\partial x_2\partial x_3}(z_1,z_2,0)
\end{split}
\end{equation}
which is the generator of the three dimensional Ornstein Uhlenbeck process reflected on the boundary of the half space $z>0$.
\end{proposition}

 \section*{Acknowledgement} 
The authors are grateful to Philippe Robert for useful conversations and exchange of notes during the preparation of this work.

  \end{document}